\documentclass[12pt]{article}

\usepackage{physor2016}

\usepackage{amsmath}
\usepackage{bm}

\newcommand{\Macro}{\ensuremath{\Sigma}}
\newcommand{\vOmega}{\ensuremath{\hat{\Omega}}}

\newcommand{\rvec}{\ensuremath{\vec{r}}}

\usepackage{graphicx}
\usepackage{tikz}
\usepgflibrary{shapes.geometric}
\usepackage{hyperref}
\hypersetup{
    colorlinks,
    citecolor=black,
    filecolor=black,
    linkcolor=black,
    urlcolor=black
}

\usepackage{booktabs}

\usepackage{siunitx}
\title{FW/CADIS-$\Omega$: AN ANGLE-INFORMED HYBRID METHOD FOR DEEP-PENETRATION RADIATION TRANSPORT}
\author{ 
  \textbf{Madicken Munk and R.~N.~Slaybaugh} \\
  Department of Nuclear Engineering, University of California, Berkeley \\
  3115B Etcheverry Hall, Berkeley, CA 94720, USA\\
  \href{mailto:madicken@berkeley.edu}{madicken@berkeley.edu}\\
  \href{mailto:slaybaugh@berkeley.edu}{slaybaugh@berkeley.edu}\\
  \\
  \textbf{Tara M.~Pandya, Seth R.~Johnson, and T.~M.~Evans}\\
  Radiation Transport Group\\
  Oak Ridge National Laboratory, P.O.\ Box 2008, Oak Ridge, TN 37831, USA\\
  \href{mailto:pandyatm@ornl.gov}{pandyatm@ornl.gov}\\
  \href{mailto:johnsonsr@ornl.gov}{johnsonsr@ornl.gov}\\
  \href{mailto:evanstm@ornl.gov}{evanstm@ornl.gov}
  }
  
\renewcommand{\shortauthor}      % Author's names here
           {M.\ Munk~et~al.}  
\renewcommand{\shorttitle}       % Short title here
           {Angle-Informed CADIS and FW-CADIS}  

\hypersetup{
  pdftitle=\shorttitle,
  pdfauthor=\shortauthor
}

\begin{document}

%\doublespacing

%\linenumbers

%------------------------------------------------------------------------------
% Make the titlepage and set the pagestyle to fancy throughout
%------------------------------------------------------------------------------
\maketitle

\begin{abstract}
A new method for generating variance reduction parameters for strongly anisotropic, deep-penetration radiation shielding studies is presented. This method generates an alternate form of the adjoint scalar flux quantity, $\phi^{\dagger}_{\Omega}$, which is used by both CADIS and FW-CADIS to generate variance reduction parameters for local and global response functions, respectively. The new method, called CADIS-$\Omega$, was implemented in the Denovo/ADVANTG software suite, and results are presented for a concrete labyrinth test problem. Results indicate that the flux generated by CADIS-$\Omega$ incorporates localized angular anisotropies in the flux effectively. CADIS-$\Omega$ outperformed CADIS in the test problem while obtaining accurate results. This initial work indicates that CADIS-$\Omega$ may be highly useful for shielding problems with strong angular anisotropies. A future test plan to fully characterize the new method is proposed, which should reveal more about the types of realistic problems for which the CADIS-$\Omega$ will be suited. 
\end{abstract}


%------------------------------------------------------------------------------
%
%------------------------------------------------------------------------------
\section{INTRODUCTION}
\label{sect::intro}

Efficiently modeling radiation transport in deep penetration shielding problems is essential to the safe operation of many types of nuclear facilities, as well as the development of monitoring and detection systems. We would like to be able to perform such calculations with Monte Carlo (MC) methods, but it can be quite challenging to obtain acceptable statistical uncertainties in the computed tallies. Thus, many variance reduction (VR) strategies, some of which we will discuss, have been developed to facilitate accurate calculations in reasonable times. 

Problems that exhibit a strong degree of angular anisotropy in particle flux tend to be even more challenging for effective VR.
Many existing VR methods do not include angular information in their techniques, and therefore do not work well for these problems.  
Some VR methods do include angular information, but, in general, these methods have narrow applicability and/or are difficult to use and/or are computationally costly\textemdash rendering them inadequate as general tools.

With the goal of developing a method to easily, reliably, and inexpensively solve fixed-source problems that exhibit a high degree of angular anisotropy, we have developed a new method that builds on the Consistent Adjoint Driven Importance Sampling (CADIS) and Forward Weighted-CADIS (FW-CADIS) methods~\cite{wagner_forward-weighted_2007}, which we will jointly refer to as FW/CADIS. 
%The new method builds on existing software infrastructure in a way that facilitates easy adoption.
FW/CADIS uses scalar flux estimates from a deterministic calculation to create VR parameters for use in MC.
Our new method uses a forward-weighted adjoint scalar flux based on a normalized contributon flux instead of a standard adjoint scalar flux, thereby including more angular information. 
We are calling the new methods CADIS-$\Omega$, FW-CADIS-$\Omega$, and FW/CADIS-$\Omega$ for both methods jointly.

We have implemented CADIS-$\Omega$ and tested it on an optically-thick source-detector labyrinth demonstration problem.
This demonstration problem has characteristics that are particularly challenging for existing VR methods and provides an opportunity to deeply investigate differences between this new method and CADIS.
We compare analog, standard CADIS, and CADIS-$\Omega$ calculations. 
Initial results show that for problems with strongly anisotropic behavior, CADIS-$\Omega$ outperforms traditional CADIS.

This paper begins with a background (Section~\ref{sect::second}) of important concepts relating to VR and existing hybrid methods for deep-penetration radiation transport. 
We then provide the context of existing methods in Section~\ref{sect::past}. 
Section~\ref{sect::methodology} describes the mathematical foundation of our proposed method and the software that we used to implement it. 
The results and accompanying discussion are described in Section~\ref{sect::results}. Section~\ref{sect::future} presents our future test plans and details of how this plan will characterize the application space beyond this single demonstration test. 
We conclude in Section~\ref{sect::conclusion}.

%------------------------------------------------------------------------------
%
%------------------------------------------------------------------------------
\section{BACKGROUND}
\label{sect::second}

Many modern radiation transport codes offer VR capabilities that employ an importance map\textemdash a measure of how important a particular region in a MC simulation is to the tally being computed\textemdash to perform VR. Generally, the importance of a cell can be defined as the ratio between the total score of particles entering the cell to the total weight of the cell \cite{booth_automatic_1982}. By generating an effective importance map for a system, the MC calculation will have faster time to solution, a reduced variance, or both. 

Importance maps can be generated a number of manners, but a particularly effective method is to use the solution, $\phi^{\dagger}$, to the adjoint formulation of the transport equation. 
The solution to the adjoint equation has long been recognized as 
% What words did you want here? I'm putting in this to have something. 
% What you wrote is good. I like it. 
directly corresponding to how influential a given source particle will be on the response function, defined as the adjoint source.
In fact, Kalos \cite{kalos_importance_1963, goertzel_monte_1958} and others have shown that an exact solution to the adjoint equation would result in a zero variance MC solution. 
The steady state, fixed-source forms of the forward and adjoint transport equation are shown in Eqs.~\ref{eq:forward} and \ref{eq:adjoint}, respectively, where $\psi(\vec{r}, \hat{\Omega}, E)$ is the angular neutron flux, $\Sigma_t$ is the total macroscopic cross section, $\Sigma_s$ is the double-differential scattering cross section, $q$ is a fixed source, and superscript $\dagger$ indicates adjoint quantities.
\begin{subequations} 
\label{eq:transport}
\begin{align}
\bigl[\hat{\Omega} \cdot \nabla + \Macro_t(\vec{r}, E)\bigr] \psi(\vec{r}, \hat{\Omega}, E)  =  \int_0^{\infty} dE' &\int_{4\pi} d\hat{\Omega'} \:\Macro_{s}(\vec{r}, E' \to E, \hat{\Omega'} \cdot \hat{\Omega}) \psi(\vec{r}, \hat{\Omega'}, E')\nonumber \\
 &+ q(\vec{r}, \vOmega, E) \label{eq:forward} \\
\bigl[-\vOmega \cdot \nabla + \Sigma_t(\rvec, E)\bigr] \psi^{\dagger}(\vec{r}, \vOmega, E) = \int_0^{\infty} dE' &\int_{4\pi} d\vOmega' \: \Sigma_s(\rvec, E \rightarrow E', \vOmega \cdot \vOmega') \psi^{\dagger}(\rvec, \vOmega', E') \nonumber \\
&+ q^{\dagger}(\vec{r}, \vOmega, E) \label{eq:adjoint}
\end{align}
\end{subequations}

%In particular, the adjoint transport equation differs from the forward equation in that particles are scattered up in energy, from E to E', and are reversed in direction, from $\vOmega$ to $\vOmega'$. 
The key to using the adjoint flux for VR is defining the adjoint source~\cite{wagner_forward-weighted_2007}. 
Physically, $\psi$ represents how the source particles go forward and affect the rest of the problem space while $\psi^{\dagger}$ represents how the particles from a source come into the solution space and affect the response. 
In this way, the angular adjoint flux represents how every part of phase space will influence the response we are seeking.
%definition of the adjoint source, $q^{\dagger}$, depends on the nature of the problem being solved. 
% changing answer to response to keep language continuous throughout the paper
For a simple source-detector problem, for example, where the desired response is some reaction rate in the detector, the adjoint source is defined as the detector response function, or $q^\dagger = \Sigma _{ d }$. 
% Adding some wording to be a little more verbose. 
%Thus, the adjoint particles start at low energy at the detector location, move away from the adjoint source (the detector location), and scatter up in energy. 

CADIS-$\Omega$ uses the contributon flux, defined in Eq.~\eqref{eq.Cont-Flux}, where contributons are pseudo-particles that carry ``response" from the radiation source to a detector ~\cite{williams_generalized_1991,williams_contributorn_1992,williams_contributon_study}. 
\begin{equation}
\Psi (\vec {r},\:\hat\Omega ,E) = \psi^{\dagger} (\vec {r},\:\hat\Omega ,E) \psi(\vec {r} ,\:\hat\Omega,E)
\label{eq.Cont-Flux} 
\end{equation}
The contributon flux includes both forward and adjoint information, expressing the importance of a particle that is born at a forward source and moves through space towards an adjoint source, contributing to the solution.
An importance map generated by a contributon flux will assign high importance to particles that are generated at the forward source and likely to generate a response in the detector. 
Becker's thesis~\cite{becker_hybrid_2009} aptly points out that this understanding of contributons is illustrated by a source-detector problem, where the forward source has little importance to the adjoint source, but does have importance to the problem solution.
In comparison, the adjoint flux places high importance on particles that are close to the adjoint source without accounting for the impact of forward transport. 

\section{PAST WORK}
\label{sect::past}

The generation of effective importance maps for automated VR techniques remains an active area of research. A wide variety of approaches have been investigated which utilize aspects of either adjoint theory or contributon theory to gain insight into importance function parameters for deep-penetration radiation shielding problems. Some techniques utilize an initial MC estimate of importance, and others utilize a determinstic calculation. In this section, we will briefly describe some of the methods utilized to generate biasing parameters for VR techniques, and consider the strengths and weaknesses of each method.

As mentioned in Section \ref{sect::second}, the likelihood of particles in a cell contributing to a tally can be calculated in a variety of manners. Booth \cite{booth_automatic_1982} showed that the importance of a cell could be estimated with MC by calculating the total score of particles entering the cell and dividing it by the total weight of the cell. Concurrently, Hendricks \cite{hendricks_code-generated_1982} developed an automated weight window generator for MCNP~\cite{brown_mcnp_2002}. The weight window generator automated VR by calculating the total weight entering and exiting the weight window target region, which allowed the generator to iteratively converge upon weight window target values. Several other methods exist to iteratively calculate weight window parameters with MC calculations, but in general, target weights in weight window maps are the inverse of the importance values for a given point in phase space.

Perhaps the most widespread and accessible of methods that utilize a deterministic transport calculation to generate an importance map are the CADIS \cite{wagner_automatic_1997,wagner_automated_1998,haghighat_monte_2003} and the FW-CADIS \cite{wagner_forward-weighted_2007,wagner_forward-weighted_2009,wagner_forward-weighted_2010} methods. These methods optimize MC transport for localized and global response functions by using consistently biased source particles and particle weights. This biasing is done using a determinstic solution for the adjoint scalar flux, $\phi^{\dagger}$, as a measure of the importance. Equations \eqref{eq:CADISmethod} and \eqref{eq:fwCADISmethod} describe the biasing parameters generated by CADIS and FW-CADIS, respectively. The biasing parameters used by CADIS are generated from the solution of the adjoint transport equation, where $\phi^{\dagger}$ is obtained from a calculation where $q^\dagger $  is set to the response function:

\begin{subequations} 
\label{eq:CADISmethod} 
\begin{equation}
\hat{q}  = \frac{\phi^{\dagger}(\vec {r} ,E)q(\vec {r} ,E)}{\iint\phi^{\dagger}(\vec {r} ,E)q(\vec {r} ,E) dE d\vec{r}} \\
         = \frac{\phi^{\dagger}(\vec {r} ,E)q(\vec {r} ,E)}{R} ,
\label{eq:weightedsource}
\end{equation}
\begin{equation}
w_0  = \frac{q}{\hat{q}} \\
     = \frac{R}{\phi^{\dagger}(\vec {r} ,E)} , 
\label{eq:startingweight}
\end{equation}
\begin{equation}
\hat{w} = \frac{R}{\phi^{\dagger}(\vec {r} ,E)}  ,
\label{eq:WW}
\end{equation}
\end{subequations}

where $\hat{q}$ is the biased source distribution, $w_0$ is the starting weight of the particles, $\hat{w}$ is the target weight of the particles, and R is the response of interest.

While CADIS is quite effective at reducing the variance in localized responses, FW-CADIS was developed to reduce the variance in global responses. The concept behind FW-CADIS is that if the calculated response can be quantified with an even number of particles sampled in each phase-space region, then the uncertainty distribution will be even across the global response tally. Rather than setting  $q^\dagger = \Sigma _{ d }$ as CADIS does for source-detector problems, FW-CADIS generates a $q^\dagger$ based on a deterministic forward calculation (Eqs. \eqref{eq:energy_space_flux}-\eqref{eq:total_doserate}). This $q^\dagger$ is used to solve the adjoint equation, and then source biasing and particle weights are adjusted with the same methods described in Eqs. \eqref{eq:weightedsource}-\eqref{eq:WW}.

\begin{subequations} 
\label{eq:fwCADISmethod}
The biasing parameters that are used by FW-CADIS depend on the desired optimized results.  
To calculate the energy and spatially dependent flux, $\phi(\vec{r},E)$, FW-CADIS sets 
\begin{equation}
q^{\dagger}  = \frac {1}{\phi(\vec{r},E)} ,
\label{eq:energy_space_flux}
\end{equation}
to calculate spatially dependent flux, $\int\phi(\vec{r},E)dE$ , 
\begin{equation}
q^{\dagger}  = \frac {1}{\int\phi(\vec{r},E)dE} , 
\label{eq:total_flux}
\end{equation}
and to calculate spatially dependent dose in the system, $\int\phi(\vec{r},E)\sigma_d(\vec{r},E)dE$ :
\begin{equation}
q^{\dagger}  = \frac {\sigma_d(\vec{r},E)}{\int\phi(\vec{r},E)\sigma_d(\vec{r},E)dE} . 
\label{eq:total_doserate}
\end{equation}
\end{subequations}

For problems with strong anisotropies in the particle flux, the importance map and biased source developed using the space/energy treatment above may not represent the real importance well enough to sufficiently accelerate the MC calculation. Note that because the scalar adjoint flux is used in Eqs.~\eqref{eq:CADISmethod} and \eqref{eq:fwCADISmethod}, the angular dependence of the importance function is not retained. The drawback of using scalar flux quantities to calculate these parameters is that no information is retained on how particles move towards the response function. However, if the angular dependence of the importance function were retained, the map would be very large (tens or hundreds of GB) and more costly to use in the MC simulation. It is worth noting that the automatic weight window generator in MCNP does have the capability of angular biasing. For some problems, this generator works well, but it requires iteration and inherently depends on a MC solution for weighting parameters, which is limited by statistical precision and may take time to acquire a weight window with an acceptable precision.

A variety of methods have incorporated angular information, to some degree, into deterministically-informed MC VR. Of note are the AVATAR \cite{van_riper_avatarautomatic_1997} and Local Importance Function Transform (LIFT) \cite{turner_automatic_1997} methods. The LIFT method uses an approximation of the adjoint flux that is piecewise continuous in angle. However, LIFT only captures linearly anisotropic scattering \cite{turner_automatic_1997-1}, which is not suitable for all types of strongly anisotropic systems. The AVATAR method, which is an implementation of the maximum entropy distribution, assumes that the flux is separable and symmetric about the average current. While AVATAR does generate biased weight windows, they are not consistently biased with the source distribution (while those in CADIS and FW-CADIS are). Neither LIFT nor AVATAR perform well in voids and low-density regions \cite{turner_automatic_1997-1}. 

More recently, Peplow et al. \cite{peplow_consistent_2012}, built off of AVATAR by consistently biasing the source distribution. In this work, Peplow et al. created an angular version of CADIS both with limited directional source biasing and without directional source biasing. These methods performed well compared to the analog, but were not sufficiently accurate in capturing angular information to scale to complex problems with strong anisotropies. 

In summary, substantive work has been performed to generate VR parameters for deep penetration shielding problems with MC. These methods generate biasing functions or weight window maps from calculated importance functions. Some methods iteratively converge on cell importances with MC, which work well if the problem is sampled with sufficient precision. Other methods utilize deterministic calculations to generate VR parameters. In this realm, both CADIS and FW-CADIS perform well for deep-penetration and relatively isotropic problems. However, CADIS and FW-CADIS struggle in problems where there are strong anisotropies in the flux. Work has been performed to generate VR parameters with angular information, but these methods are either not automated, have limited accessibility, have limited applications, or struggle to capture the anisotropy for a large problem subset.

%------------------------------------------------------------------------------
%
%------------------------------------------------------------------------------
\section{METHODOLOGY}
\label{sect::methodology}

We have seen that past methods are challenged by problems with strong angular anisotropies or have limited use cases.
We build on past methods but calcualte the adjoint scalar flux in a way that has not been done before.
Using more angular information should improve performance and increase reliability for problems in which angular information is particularly important. 
In this section we describe the theory behind the new method and why we believe it will work for highly anisotropic problems. We also present an overview of the software implementation.  
% Someday work on making this section flow more smoothly. 

%------------------------------------------------------------------------------
%
%------------------------------------------------------------------------------
\subsection{Theory}
\label{subsect::theory}

Our new automated hybrid method incorporates angular information into the biasing parameters for FW/CADIS while not explicitly biasing in angle. 
That is, we generate space- and energy-dependent importance maps that incorporate the flux anisotropy in a more effective way than current implementations without adding the complication of angular weight windows. 
Both CADIS and FW-CADIS use an adjoint scalar flux to generate biasing parameters for VR as shown in Eqs.~\eqref{eq:CADISmethod} and \eqref{eq:fwCADISmethod}. 
FW/CADIS-$\Omega$ uses the same equations for weight windows and source biasing, but our method generates the adjoint flux differently. 
We first take the product of the adjoint angular flux with the forward angular flux (this quantity is the contributon flux), integrate over angle, and divide by the integrated forward angular flux as shown in Eq.~\eqref{eq:angularhybrid}.
This quantity, which we designate $\phi^{\dagger}_{\Omega}$, is then used in the FW/CADIS methods.
\begin{equation} 
\phi^{\dagger}_{\Omega}(\vec{r},E) = \frac{\int_{4\pi} \psi(\vec {r} ,E,\hat{\Omega})\psi^{\dagger}(\vec {r} ,E,\hat{\Omega})d\hat\Omega }{\int_{4\pi}\psi(\vec {r} ,E,\hat{\Omega})d\hat\Omega}
\label{eq:angularhybrid}
\end{equation}

In a strongly anisotropic system, the adjoint scalar flux generated by Eq.~\eqref{eq:angularhybrid} will be informed by what directions were most prominent in the forward case. 
We can see this by considering the contributon flux, the numerator of Eq.~\eqref{eq:angularhybrid}.
The contributon flux expresses preferential transport paths for particles that are born at the forward source and generate a response.
Therefore, the particles in $\phi^{\dagger}_{\Omega}$ include the impact from how the direction they are moving influences the answer. 
This should allow for more effective MC transport when angular effects are important. 
Note that in an isotropic system, $\phi^{\dagger}_{\Omega}$ will be essentially the same as $\phi^{\dagger}$. 

%We are interested in the contributon flux because it finds preferential transport paths for particles that are born at the forward source and generate a response. This will have a strong effect on our variance reduction parameters, as the transport paths generated by the contributon flux will also allow for preferential transport paths in our adjusted Monte Carlo simulations. This will likely lead to a faster time to the solution, as less transport sampling will be performed in regions near the detector that are unlikely to have a source of forward source particles. 

%An importance map generated by the traditional adjoint scalar flux will differ significantly from a flux generated with a contributon flux. 
%An importance map 
% Madicken -- I'm not really sure why you commented out everything here but left these three words? Do you want me to add to it? 
%from the adjoint flux will merely have greater importance with increasing proximity to the adjoint source, while the imporance map generated from 
%that includes the contributon flux will have strong importances in regions important to both the forward and adjoint angular fluxes. 
%
% 
%To calculate the quantity in eq. \eqref{eq:angularhybrid}, full angular flux maps of both the forward and adjoint problems will be required. 
%Further, an extra deterministic calculation (a forward calculation) will be required for CADIS-$\Omega$ calculations.

%------------------------------------------------------------------------------
%
%------------------------------------------------------------------------------

\subsection{Implementation}
\label{subsect::implementation}

We implemented the new method through the AutomateD VAriaNce reducTion Generator\\ (ADVANTG)~\cite{wagner_automated_2002, mosher_new_2010} software developed at Oak Ridge National Laboratory (ORNL). 
ADVANTG automates the generation of the importance map and biased source distribution created using either the CADIS or FW-CADIS methods for use in MCNP5~\cite{brown_mcnp_2002}. 
An input file in MCNP syntax is provided by the user in addition to some instructions for running ADVANTG. 
ADVANTG uses this information to exectue the discrete ordinates solver Denovo~\cite{evans_denovo:_2010} and uses the output to create the VR parameters for MCNP.
%The deterministic calculations can be performed using multiple cores and/or processors (e.g., on multi-core desktop systems and clusters). 
%ADVANTG takes Denovo's output, executes the CADIS or FW-CADIS methods, and the final variance reduction parameters are output in a format that can be used with unmodified versions of MCNP. 
%The primary objective of the development of ADVANTG has been to reduce both the user effort and the computational time required to obtain accurate and precise tally estimates across a broad range of challenging transport application areas.
The main reason we chose this system is the implementation is nearly invisible to the user and therefore the experience of using this method is nearly identical to using CADIS or FW-CADIS.
This facilitates easy adoption of software that is already easy to use.

%That is, the user simply adds an additional instruction asking to use the angle informed method and the interface does not change otherwise.
%Further, only one MCNP input file is required to compare the new method to FW/CADIS.
%Finally, it was simpler to implement the new method through ADVANTG than starting separately as we could take advantage of so much existing infrastructure in the coupling.

The major modifications required to implement this method were made to Denovo. 
The angular flux is typically not stored or written by deterministic solvers as the desired output is often the scalar flux.
The new method, however, requires both the forward and adjoint angular fluxes to create the scalar flux. 
Denovo was therefore modified to store and write the angular flux.
A new function was also added that performs the integration in Eq.~\eqref{eq:angularhybrid}. 
This set of augmented scalar fluxes is then written the same way as any scalar flux output from Denovo.

The benefit of the bulk of the implementation being in Denovo is several fold. 
From the standpoint of ADVANTG, there are very few differences between the new method and FW/CADIS, making implementation straightforward.
Further, anyone who finds a use for a scalar flux created as in Eq.~\eqref{eq:angularhybrid} will now be able  access it.
Finally, it might be useful to have access to the full angular flux. 
Examining the angular flux for a problem could have research or pedagogical implications, and a VR method that uses the angular flux explicitly could be more easily developed in the future.

%------------------------------------------------------------------------------
%
%------------------------------------------------------------------------------
\section{RESULTS AND DISCUSSION} 
\label{sect::results}

% also address the ray effects differences between our method and CADIS
% change CADIS dots in response plots 
% Consider the detector response function for different energy neutrons. 
% Could put the detector response function on a dual axes for the response plot. 
% Talk about how well they perform relative to each other. 
%% This is what we did, this is what we saw, this is how we explain it. 

To begin to characterize the performance of CADIS-$\Omega$ we chose to study a Labyrinth source-detector problem.
Such labyrinth problems contain some angle dependence that is not well captured by standard CADIS~\cite{peplow_consistent_2012}. 
The motivation for this initial test is to investigate the differences between CADIS and CADIS-$\Omega$.
We also performed an analog calculation, comparing  the tally result, the relative uncertainty, and the MC figure of merit (FOM) of the response tally of the methods.
The FOM is given below, where we used the average relative error of the tally total.
\[\text{FOM} = \frac{1}{(\text{calculation time})*(\text{average relative error})^2}\:. \]

All MCNP calculations used 10,000,000 particles, an f4 tally for the total reaction rate in the NaI detector (using a path length estimator), and energy bins matching the bounds of the library used by Denovo. 
Denovo used the 27g19n shielding cross section library (a library included in the ADVANTG software distribution), a step characteristic (SC) spatial solver, quadruple range (QR) quadrature with four azimuthal and four polar angles per octant, a $P_N$ order of 3, and 135,000 spatial voxels.% In ADVANTG, the quadruple range quadrature is set with four azimuthal and four polar angles per octant, and a default quadrature order of 10.

\begin{figure}
  \begin{center}
    \includegraphics[width=0.80\textwidth]{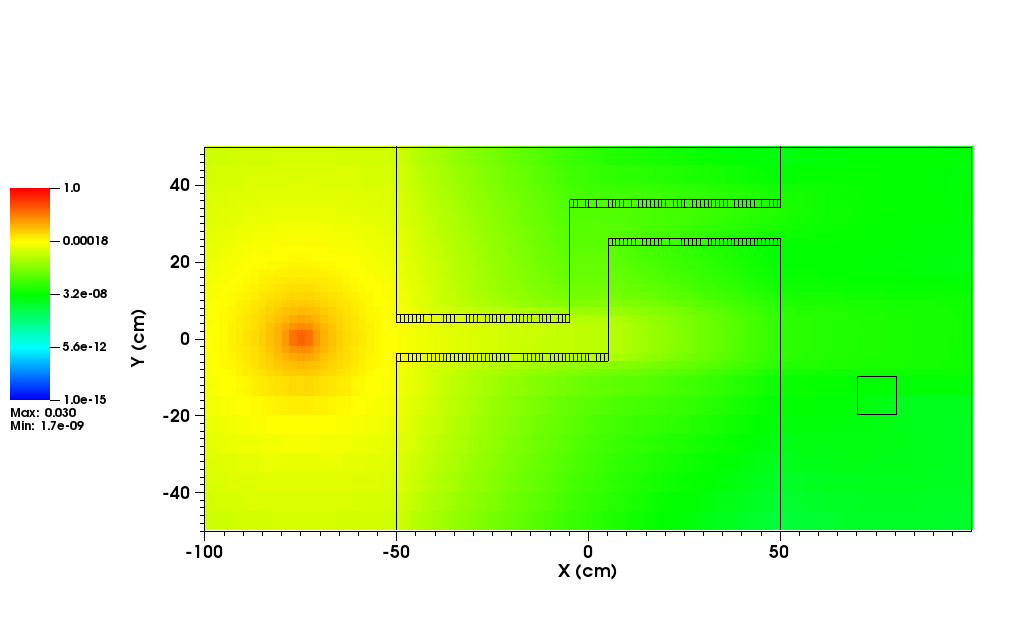}
    \caption[]{\label{fig::fwdflux}Forward deterministic flux generated by Denovo for Labyrinth problem.}
  \end{center}
\end{figure}
% the axes and labels are really hard to read 
Figure~\ref{fig::fwdflux} shows the problem geometry overlaid on a forward deterministic flux map. The labyrinth is comprised of a concrete shield, a 10 MeV isotropic point source at (-75, 0, 50), and a 10x10x10 cm cubic NaI detector on the other side of the shield centered at (75, -15, 50). The edges of the problem have vacuum boundary conditions, and the cross section libraries for MCNP were continuous energy. 

The deterministic portion of the CADIS method took 41.5 minutes and CADIS-$\Omega$ took 83.0 minutes, both on a single core of a 3 GHz Intel Core i7 processor. The MC runtimes for the analog, CADIS, and CADIS-$\Omega$ were 64, 483.4, and 408.9 minutes, respectively, on a single, 20-core node of Intel Xeon E5-2670, resulting in FOMs of 493, 5.1, and 145.0, respectively. These results are summarized in Table \ref{tab:FOMLabI}. In observing the adjusted FOM values, one should note that the system that performed the deterministic calculations was different than that of the MC calculations. 

\begin{figure}
  \begin{center}
    \includegraphics[width=0.49\textwidth]{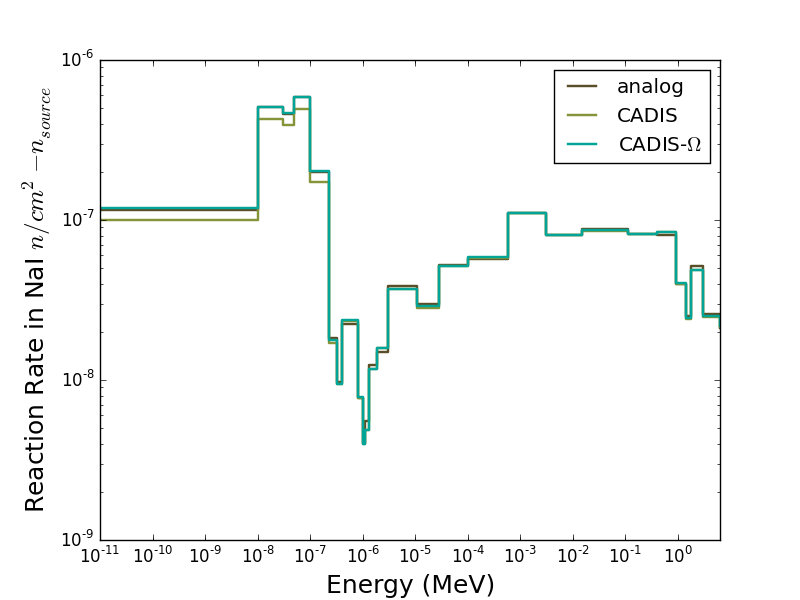}
    \includegraphics[width=0.49\textwidth]{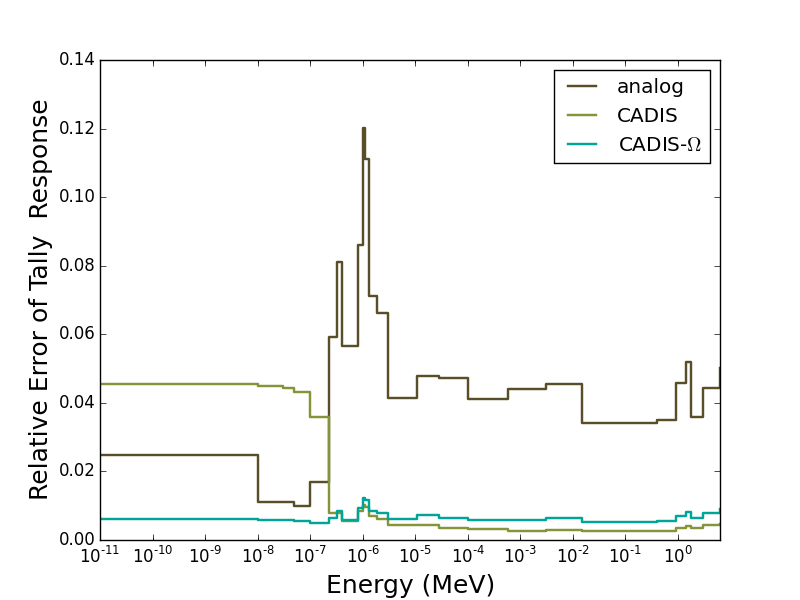}
    \caption[]{\label{fig::tallyresponse} Responses (left) and relative uncertainties (right) in NaI detector tally for the Labyrinth problem. }
  \end{center}
\end{figure}
% Axes aren't great, hard to read the CADIS line
% These need error bars if possible, otherwise it is very difficult to compare

\begin{figure}
  \begin{center}
    \includegraphics[width=0.49\textwidth]{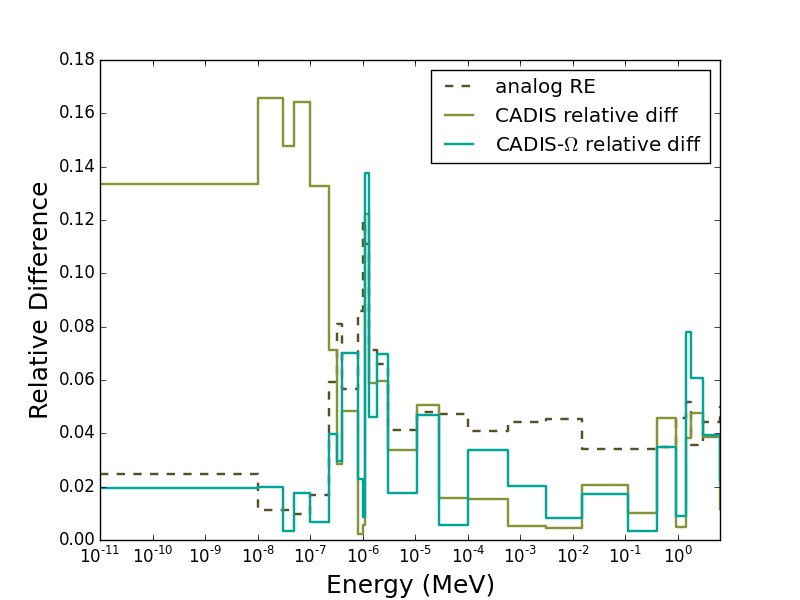}
    \caption[]{\label{fig::rediffs} Relative differences between the CADIS and CADIS-$\Omega$ responses and the analog response. For scale, the analog relative error is also included. }
  \end{center}
\end{figure}

The tally responses and relative uncertainties for the NaI detector in the Labyrinth problem are shown in Fig.~\ref{fig::tallyresponse} for the analog, CADIS, and CADIS-$\Omega$ methods.  CADIS-$\Omega$ has a lower relative error than CADIS for lower energy bins, but a slightly higher relative error for higher energy bins. The relative error for CADIS-$\Omega$ appears to be uniformly low. Comparatively, CADIS has the highest relative errors of all three methods at low neutron energies, but the lowest relative error for the highest energies. Not unsurprisingly, the analog response had the highest uncertainty in the energy bins with the lowest response. Adding error bars on the response plot of Figure~\ref{fig::tallyresponse} obscured the data. In general, the responses are equivalent to within one standard deviation. This is not the case at low energies, where none of the methods lie within one standard deviation of each other, and CADIS lies well outside of the bounds of both the analog and CADIS-$\Omega$ methods. 
To facilitate comparison, the relative difference between the CADIS and CADIS-$\Omega$ responses and the analog response is shown in Fig.~\ref{fig::rediffs}. 

\begin{figure}
  \begin{center}
    \includegraphics[width=0.49\textwidth]{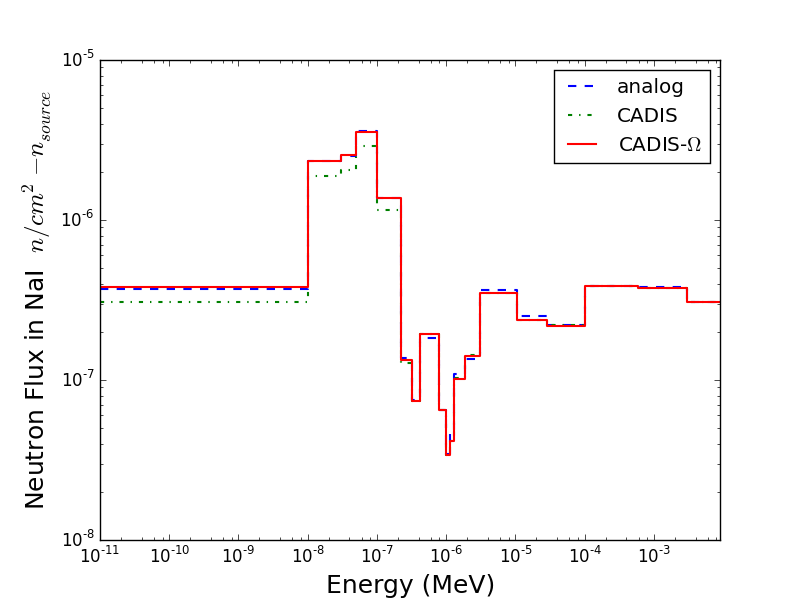}
    \includegraphics[width=0.49\textwidth]{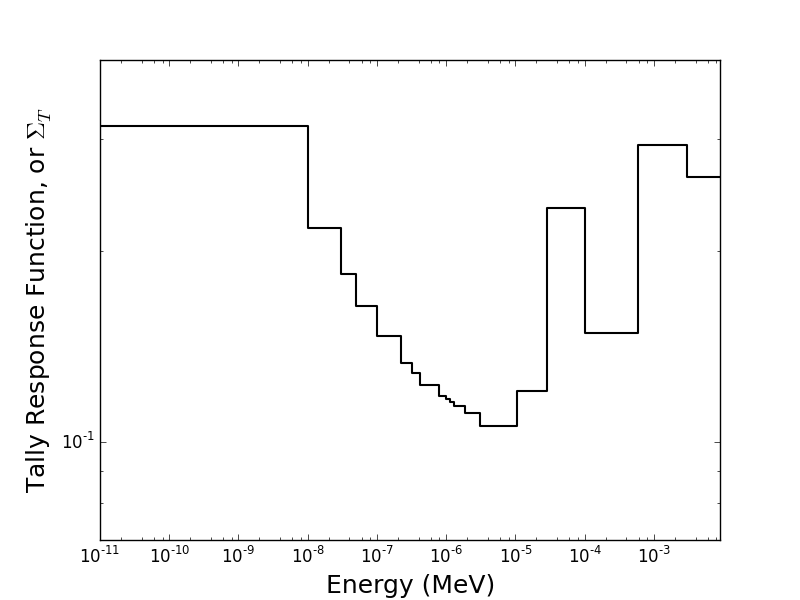}
    \caption[]{\label{fig::tallyproducts} Neutron flux at the NaI detector (left) and the sampled response function (right) that are folded together to obtain detector response. }
  \end{center}
\end{figure}

By separating the detector response into the scalar flux and the response function ($\Sigma_{T}$), shown in Fig.~\ref{fig::tallyproducts}, we gain some insight as to why the results in the low energy region are in disagreement.
The issues that lead to both CADIS and CADIS-$\Omega$'s disagreement with the analog calculation at low energies may be that the behavior of the intermediate energy range neutrons is not being effectively captured in the VR because there are so few of them and they do not dramatically impact the detector response. 
The source of this discrepancy merits further investigation. 
% Not sure if this rearrangement is better, but I think it puts the discussions closer to the figures?

\begin{figure}
  \begin{center}
    \includegraphics[width=0.80\textwidth]{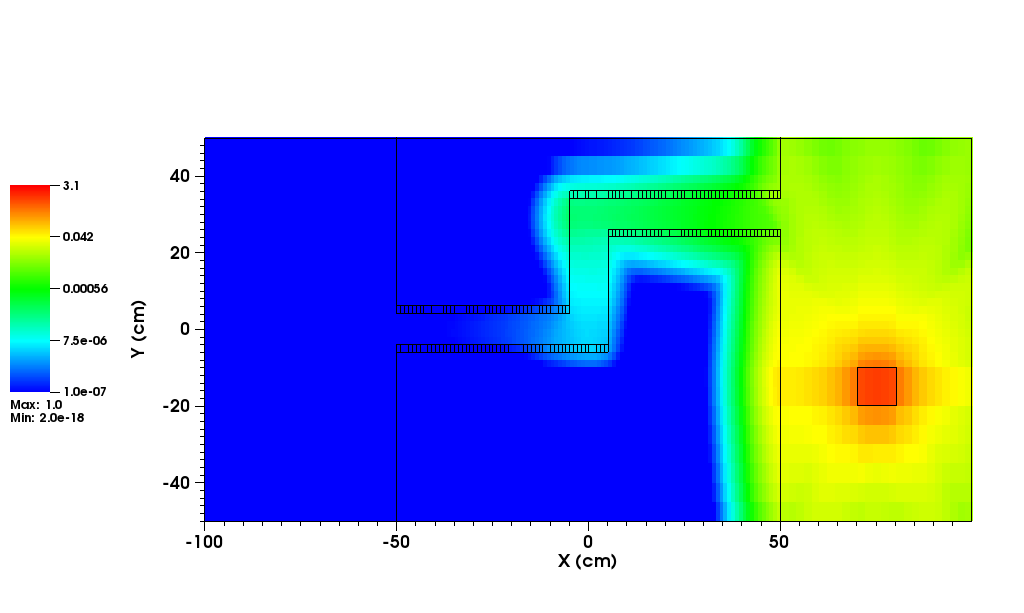}
    \includegraphics[width=0.80\textwidth]{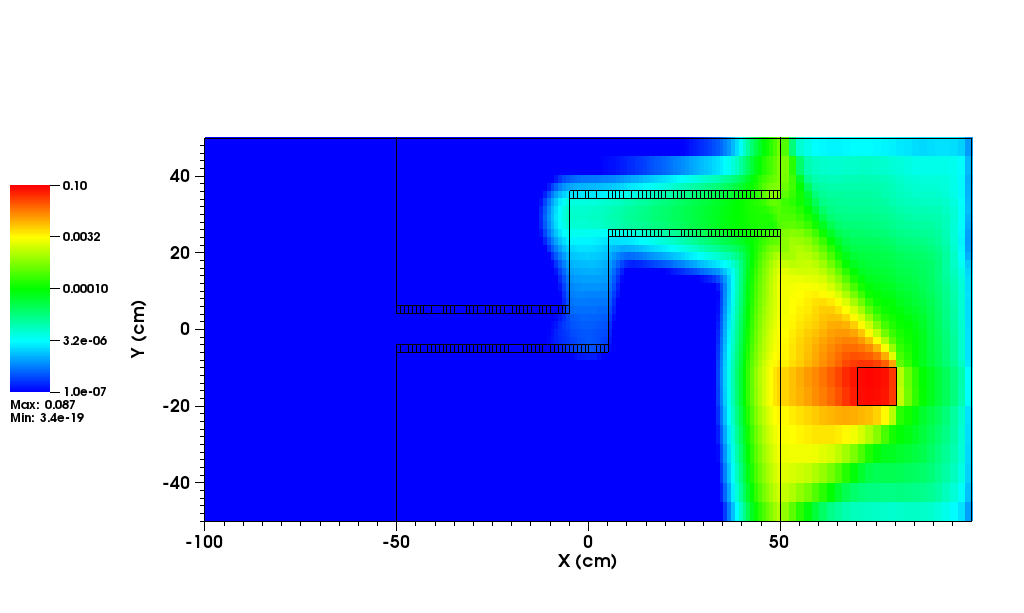}
    \caption[]{\label{fig::adjoint_fluxes_group26} Group 26 adjoint flux generated by standard Denovo (top) and the CADIS-$\Omega$ method (bottom).}
  \end{center}
\end{figure}

\begin{figure}
  \begin{center}
    \includegraphics[width=0.80\textwidth]{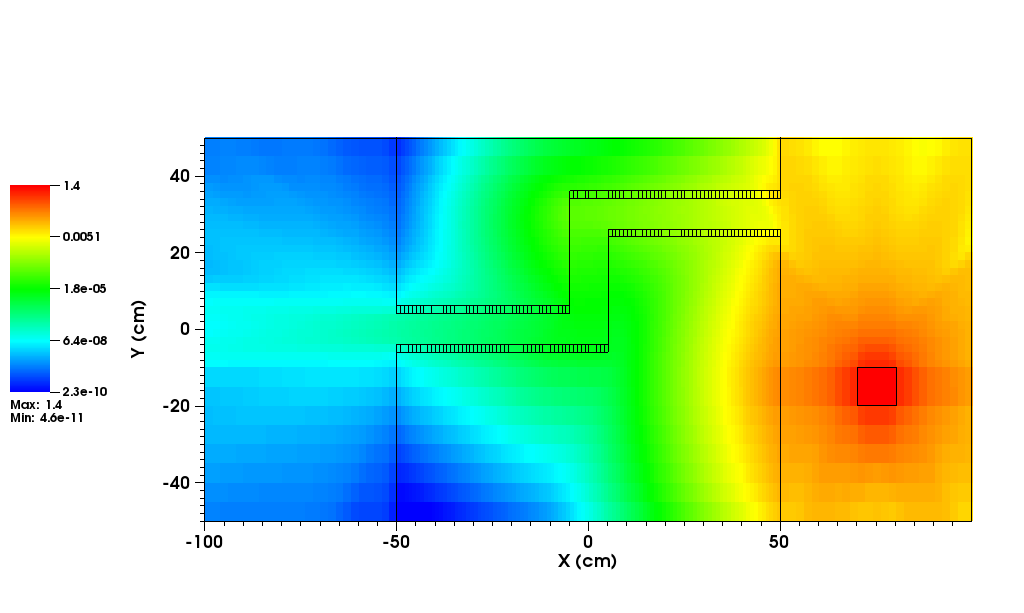}
    \includegraphics[width=0.80\textwidth]{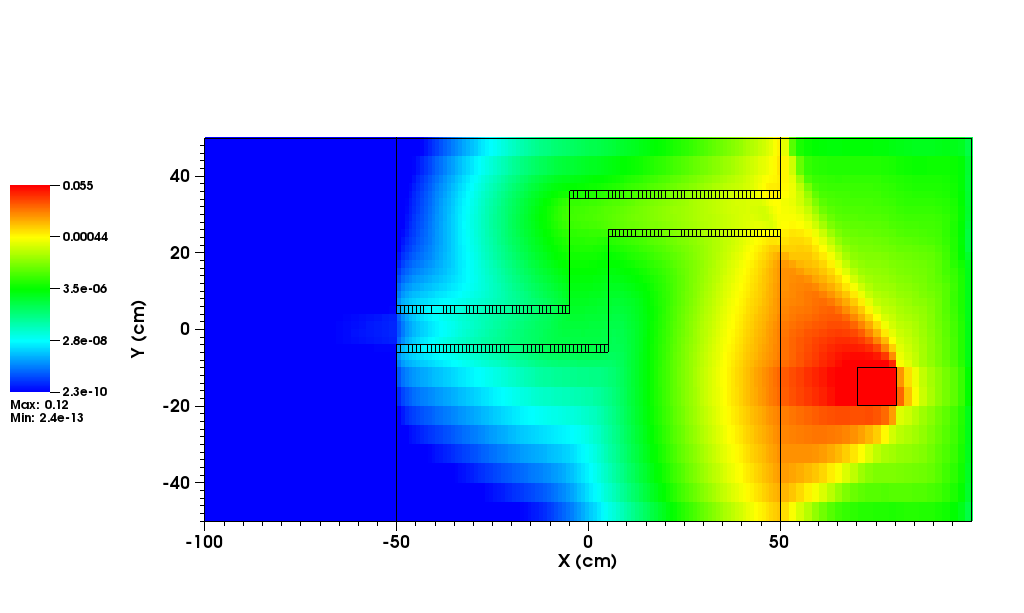}
    \caption[]{\label{fig::adjoint_fluxes_group11} Group 11 adjoint flux generated by standard Denovo (top) and the CADIS-$\Omega$ method (bottom).}
  \end{center}
\end{figure}

To investigate the reason that CADIS-$\Omega$ has higher relative error in some energies, we compared the adjoint flux distributions generated by CADIS and CADIS-$\Omega$ in a lower energy group, group 26 (1E-11 to 1E-08 MeV), seen in Fig.~\ref{fig::adjoint_fluxes_group26}, and a higher energy group, group 11 (2.9E-05 to 1.01E-04 MeV), seen in Fig.~\ref{fig::adjoint_fluxes_group11}.
Note that in these figures the color scales are not the same between groups so that the relative behavior of the flux within a method can be compared between the methods. It is the absolute difference between the flux values that will influence the weight window bounds, so this type of approach is reasonable. 

In group 26, CADIS-$\Omega$ creates an adjoint flux in which the majority of the group 26 neutrons that contribute to the detector response are those that exit the channel in the small angle between the shield and the detector, while the standard adjoint indicates fairly symmetric contribution. 
% are reflected off of the back and bottom walls, 
% while the standard adjoint indicates  fairly symmetric contribution. 
This illustrates that CADIS-$\Omega$ captures angular information differently than CADIS. 
We expect the shape of the CADIS map given the non-inclusion of angular information. 
The CADIS-$\Omega$ map seems to be strongly dominated by the streaming behavior of the particles between the shield and the detector in the right half of the problem, indicating a locally anisotropic behavior in this region. This is likely due to the lack of scattering reactions off of air for neutrons in this region. As particles travel beyond the detector, their probability for inducing a response decreases as their probability for exiting the problem increases. The decrease in particle importance near the problem boundaries is consistent with what has been observed for the contributon flux. Further, the importance of low energy neutrons contributing to the detector response drops off sharply shortly into the shield. This is likely due to the high probability of absorption of low energy neutrons into the shield. That is, low energy neutrons that are in the shield have a much lower survival probability than higher energy neutrons, like those in group 11. 
As a result, the majority of the group 26 neutrons that stream directly towards the detector once exiting the channel will contribute to the response, resulting in a locally anisotropic response-generating flux. 

This asymmetric effect in the CADIS-$\Omega$ case exists, but is not as prominent in group 11. 
Unlike group 26, the group 11 neutrons are far more likely to survive an interaction with the shield. As a result, the response-generating flux of group 11 has a similar probability to be either from streaming out of the channel or from exiting the shield, resulting in a more locally-isotropic flux near the detector volume. As a result, a much greater portion of the shield contains particles of high importance to the detector response. Therefore, the group 11 CADIS-$\Omega$ adjoint flux is less dominated by the streaming behavior of particles exiting the channel. 
% Great description in both paragraphs!

 \begin{table}
  \centering
  \caption{\label{tab:FOMLabI}FOMs and Timings for Labyrinth Problem}
  \begin{tabular}{l|cc|cc}
    \toprule
        & FOM$_{MC}$ & FOM$_{adjusted}$ & T$_{MC}$ (minutes) & T$_{determinstic}$ (minutes) \\
    \hline
    Analog           & 493   &  493     & 64.9      & 0.00 \\ 
    CADIS            & 5.1   &  ---     & 483.4     & 41.5  \\
    CADIS-$\Omega$   & 145.0 &  -----   & 408.9     & 83.0  \\  
	\bottomrule
  \end{tabular}
\end{table}

As mentioned previously, Table \ref{tab:FOMLabI} reports FOMs using the Monte Carlo simulation time alone (designated MC), as well as an adjusted FOM that incorporates the runtimes for deterministic transport (designated adjusted). For a simple problem with the refined angular and spatial discretizations we chose, the deterministic calculation times were comparatively significant.
Because CADIS-$\Omega$ uses two deterministic transport calculations while CADIS uses one, incorporating the increased time for the calculation is necessary for a fair comparison.
It is worth noting that FW-CADIS-$\Omega$ and FW-CADIS would no longer have this discrepancy in deterministic calculation time.

The FOMs reported in Table \ref{tab:FOMLabI} reveal the significant effect that the high relative errors on low energy bin responses had on the FOM for CADIS. Because the low energy bins had a high response cross section in addition to a high flux, these bins strongly affected the FOM for CADIS. The large difference in flux anisotropy near the detector in Fig.~\ref{fig::adjoint_fluxes_group26} may also explain this. The biasing parameters generated by CADIS are encouraging the MC calculation to spend equal calculation power splitting and rouletting for all directions entering the detector even though the flux is most likely to stream directly from the channel to generate a response. As a result, CADIS spends more time tracking in regions that do not contribute to the solution. Conversely, CADIS-$\Omega$ does capture this behavior, suggesting it will perform well for problems with strong, local angular anisotropy.  

% We need to mention the comparison to analog. Here's a shot at expalaining; tweak as necessary. 
It is noteworthy that, for this problem, analog MC produced a higher FOM than either VR method. 
It may be that this problem was not sufficiently challenging to see the benefit of VR compared to analog. 
It is, however, important to mention that even though the FOM for analog is highest, the results have unacceptably high relative error in a fair amount of the energy spectrum. 
It may be that we would rather have the uniformly low response from CADIS-$\Omega$ at a slightly higher time cost.
Running more analog particles to see how long it takes to reduce the relative error to acceptable levels everywhere will be an informative task. 

Another important concept that will affect the success of both CADIS and CADIS-$\Omega$ are ray effects, which are manifestations of rays in the flux that result from discretizing the transport calculation in angle. Both Figs.~\ref{fig::adjoint_fluxes_group26} and ~\ref{fig::adjoint_fluxes_group11} exhibit some ray effects behind the shield. Comparing the flux maps between CADIS and  CADIS-$\Omega$, CADIS-$\Omega$ does not seem to suffer any more from ray effects than traditional CADIS and even appears to reduce them in some locations. This softening of the rays may be caused by the integration performed in the numerator of Eq.~\eqref{eq:angularhybrid}. In this problem, the streaming path out of the labyrinth is perpendicular to the rays generated by the adjoint source, resulting in a softened adjoint-$\Omega$ flux map. However, a future problem where the ray effects in both the forward and adjoint sources lie parallel to one another might result in synergistic ray effects, which would have adverse effects on the VR effectiveness of the $\Omega$-methods. The extent to which CADIS-$\Omega$ can reduce ray effects is an important attribute to study in future tests.

Overall, the performance of CADIS-$\Omega$ compared to CADIS and analog MC is promising. In our simple demonstration test problem, CADIS-$\Omega$ performed significantly better than CADIS by both FOM and relative error distribution inspection. CADIS-$\Omega$ had a more uniform uncertainty distribution in the response tally than both CADIS and the analog calculation. To determine the cause of these results, we compared the adjoint flux quantities calculated by CADIS and CADIS-$\Omega$. Inspection of the flux maps revealed that CADIS-$\Omega$ captures anisotropies in the flux, and these anisotropies are reflected in the VR parameters that are generated. Based on both theory and the method's performance in this demonstration problem, we expect that this method will perform well in problems with strong flux anisotropy generated by streaming paths or material heterogeneity, while performing poorly in problems with significant scattering effects where the flux may become more locally isotropic.

%------------------------------------------------------------------------------
%
%------------------------------------------------------------------------------
\section{FUTURE WORK} 
\label{sect::future}
 
Section \ref{sect::results} showed that CADIS-$\Omega$ can outperform CADIS in a problem with a strong degree of angular anisotropy in the flux. However, there are a number of physical means by which angular anisotropy in the flux might occur. To further characterize the performance of CADIS-$\Omega$ and FW-CADIS-$\Omega$ more fully, we are planning to use a suite of test problems that cover a wide range of physical characteristics. The testing scheme identified in Table \ref{tab:testprobs} summarizes our testing plans. As noted in Section \ref{sect::results}, problems where the direction of particles is important to the result are where we anticipate the strength to lie in CADIS-$\Omega$. Our test problems aim to address all of the physical means by which the flux could become strongly anisotropic, and some have shown to be difficult for traditional FW/CADIS~\cite{peplow_consistent_2012}. FW/CADIS-$\Omega$'s sensitivity to various determinstic discretization parameters will also be investigated. 

 \begin{table}
  \centering
  \caption{Proposed Test Problem Coverage}
  \begin{tabular}{l|cccc}
    \toprule
    Problem Name & \multicolumn{4}{c}{Problem Coverage} \\
    \hline
     & Streaming Paths & High Scatter & Highly Heterog. & Beam Problem \\
    \hline
    Streaming Channel   & X & & & X \\ 
    Metal Plate         & X & X & X &  \\
    Labyrinth Variants  & X & X & X &  \\ 
    Spherical Boat      & X & & X & X \\  
    Kobayashi Benchmark & X & X &  &  \\   
	\bottomrule
  \end{tabular}
  \label{tab:testprobs}
\end{table}

After FW/CADIS-$\Omega$ has been fully characterized, it will be tested with large problems of interest. 
Specifically, we are interested in real problems that contain strong angular anisotropies and that are difficult to solve, such as active interrogation problems, used fuel cask storage (through air vents and between casks on ISFISIs), and facility calculations containing long gaps or streaming materials embedded in moderators. 
These large problem studies will be informed by the small problem studies and demonstrate the potential impact of this new method. 

%------------------------------------------------------------------------------
%
%------------------------------------------------------------------------------
\section{CONCLUSION} 
\label{sect::conclusion}

The existing fleet of fully automated hybrid methods does not include a method that has the ability to effectively capture angular anisotropies in the flux for effective VR in deep-penetration radiation shielding problems. For problems that are strongly anisotropic, existing automated hybrid methods have poor performance, and do not produce results with the quality necessary to be used on a large scale. There has been a strong push to remedy this in the hybrid methods community, but the methods that incorporate angular information for variance reduction are either not automated, not widely accessible, or are limited to particular anisotropic problem types because of assumptions in their methodology. 

We have presented a new method, FW/CADIS-$\Omega$, that includes information from the full energy- space- and angle-dependent angular fluxes from a deterministic calculation to generate a weighted form of the adjoint scalar flux. This weighted form of the adjoint flux is calculated by postprocessing the adjoint angular flux values from the discrete ordinates solver Denovo. This scalar flux value, $\phi^{\dagger}_{\Omega}$, is then passed to ADVANTG to use in CADIS and FW-CADIS. By using this scalar flux quantity, the weight windows generated by FW/CADIS-$\Omega$ will incorporate angular information without explicitly using angular biasing techniques. This method is based upon a strong theoretical foundation; it incorporates concepts from both adjoint and contributon theory. Because it has been implemented in production-level software, it will be easy to use and adopt. 

We have explored the success of CADIS-$\Omega$ in a simple Labyrinth test problem and found promising results. 
Compared to CADIS, CADIS-$\Omega$ had a significantly higher FOM, a more uniform uncertainty distribution, and agreed well with analog results. 
In comparing $\phi^{\dagger}$ and $\phi^{\dagger}_{\Omega}$ in two energy groups, we found that CADIS-$\Omega$ generated dramatically different adjoint scalar flux values in regions with strong local anisotropy. 
This initial test showed that the $\phi^{\dagger}_{\Omega}$ captures local anisotropies in the flux, which is reflected in the variance reduction parameters generated by CADIS-$\Omega$.
%This is one of the contributing factors that explains CADIS-$\Omega$'s stronger performance for the example problem. 
We view these initial results as a good indication that FW/CADIS-$\Omega$ will be an effective automated variance reduction tool for problems with strong angular anisotropies and that it may significantly improve our ability to solve difficult real-world shielding problems with these characteristics. 

%While the example problem begins to demonstrate the difference between CADIS and CADIS-$\Omega$, we have developed a comprehensive testing plan that will help to characterize FW/CADIS-$\Omega$ more fully. The testing phase of this project aims to (1) characterize the method's performance in a variety of problems that have anisotropies, and (2) demonstrate the method's applicability for large, realistic problems. 
%
%  Furthermore, the potential impact of creating a successful effective automated variance reduction parameters for problems with strong angular anisotropies in the flux are numerous. Large, site-level simulations with localized flux anisotropies will be solvable. Small, active-interrogation response functions can be quickly estimated at inspection stations. Shielding design and optimization will be faster. By  presenting, developing, characterizing, and testing FW/CADIS-$\Omega$, we are improving the understanding of how to incorporate angular information into variance reduction parameters. 

%I'm not super happy with that last sentence, but I wasn't sure where to take it. I'll keep thinking about how to change it. 

%------------------------------------------------------------------------------
%
%------------------------------------------------------------------------------
\section*{ACKNOWLEDGMENTS}

This material is based on work supported by the Department of Energy under award number DE-NE0008286. This report was prepared as an account of work sponsored by an agency of the United States Government. Neither the United States Government nor any agency thereof, nor any of their employees, makes any warranty, express or implied, or assumes any legal liability or responsibility for the accuracy, completeness, or usefulness of any information, apparatus, product, or process disclosed, or represents that its use would not infringe privately owned rights. Reference herein to any specific commercial product, process, or service by trade name, trademark, manufacturer, or otherwise does not necessarily constitute or imply its endorsement, recommendation, or favoring by the United States Government or any agency thereof. The views and opinions of the authors expressed herein do not necessarily state or reflect those of the United States Government or any agency thereof.

\bibliographystyle{physor2016}
\bibliography{myphysor2016}

\appendix

\makeatletter
\def\@seccntformat#1{APPENDIX \csname the#1\endcsname.~}
\makeatother

%------------------------------------------------------------------------------
% If you need to make one (or more) appendix (appendices), place them here as
% sections
%%------------------------------------------------------------------------------
%\section{HOW TO MAKE APPENDICES}
%\label{app::a}
%
%This is a placeholder for my first appendix
%
%\section{OTHER APPENDIX STUFF}
%\label{app::b}
%
%This is a placeholder for my second appendix

\end{document}